\documentclass[11pt]{article}
\usepackage[usenames]{color}
\usepackage{amssymb}
\usepackage[colorlinks=true,
linkcolor=webgreen,
filecolor=webbrown,
citecolor=webgreen]{hyperref}

\definecolor{webgreen}{rgb}{0,.5,0}
\definecolor{webbrown}{rgb}{.6,0,0}

\newtheorem{theorem}{Theorem}

\newcommand{\bsq}{{\vrule height .9ex width .8ex depth -.1ex }}
\newcommand{\eeq}{\end{equation}}
\newcommand{\beql}[1]{\begin{equation}\label{#1}}
\newcommand{\eqn}[1]{(\ref{#1})}

\begin{document}
\begin{center}

{\Large{\bf The Number of Hierarchical Orderings}} \\
\vspace{1\baselineskip}
{\em N. J. A. Sloane} \\
AT\&T Shannon Labs \\
 Florham Park, NJ 07932--0971, USA \\
(njas@research.att.com, http://www.research.att.com/$\sim$njas) \\
\vspace{1\baselineskip}
{\em Thomas Wieder} \\
Darmstadt University of Technology \\
Institute for Materials Science \\
D--64287 Darmstadt, Germany \\
(thomas.wieder@epost.de, http://homepages.tu-darmstadt.de/$\sim$wieder) \\
\vspace{2\baselineskip}
{\bf Abstract} \\
\vspace{.4\baselineskip}
\end{center}

An ordered set-partition (or preferential arrangement) 
of $n$ labeled elements represents a single ``hierarchy'';
these are enumerated by the ordered Bell numbers.
In this note we determine the number of ``hierarchical orderings''
or ``societies'',
where the $n$ elements are first partitioned into $m \le n$ subsets
and a hierarchy is specified for each subset.
We also consider the unlabeled case, where the ordered Bell numbers
are replaced by the composition numbers.
If there is only a single hierarchy, we show that the average
rank of an element is asymptotic to $n/(4 \log 2)$ in the labeled case
and to $n/4$ in the unlabeled case.

\begin{center}{Keywords: ordered set-partition, preferential arrangement, hierarchical ordering, society}\end{center}
\begin{center}{AMS 2000 Classification: Primary 06A07}\end{center}
\vspace{.4\baselineskip}

\section{Introduction}
Suppose we are given a set $S$ of $n$ labeled elements (or ``individuals'').
The number of ordered set-partitions (or ``hierarchies'') on $S$
is given by the ordered Bell number $B_n$
(see sequence \htmladdnormallink{A670}{http://www.research.att.com/cgi-bin/access.cgi/as/njas/sequences/eisA.cgi?Anum=A000670} in \cite{OEIS} for properties and references).
A {\em hierarchical ordering} or {\em society}
on $S$ is specified by
first distributing the elements into $m \le n$ unlabeled and
nonempty subsets, and forming an ordered set-partition
in each subset.
In Section \ref{Sec_L} we will determine the number of 
different structures of this type.
Section \ref{Sec_U} discusses the analogous question
in the unlabeled case, when the
elements are indistinguishable.

The original motivation for this work was to try to describe the structure
of a ``typical'' society. Of course Sections \ref{Sec_L} and \ref{Sec_U}
just enumerate them. However, in Section \ref{Sec_H}
we consider the distribution of ranks in a random
selection of a single hierarchy, in both the
labeled and unlabeled cases.

Structures of the kind considered here were discussed 
in a classic combinatorics paper by Motzkin \cite{Motz71}.
In his terminology, the labeled structures in Section \ref{Sec_L}
would be called ``sets of lists of sets'', and the unlabeled
structures in Section \ref{Sec_U} ``sets of lists of numbers''.\footnote{The
Maple package {\em combstruct}
\cite{INRIA}
makes it easy to generate such structures.
For example, the labeled structures in Section \ref{Sec_L}
are defined by the specification
$$
[H, \{H = Set(Sequence(Set(Z, card >= 1), card >= 1))\}, labeled].
$$
\cite{OEIS} now contains the sequences enumerating all
such Motzkin-type structures
in which the specification involves up to three
occurrences of {\em Set} and {\em Sequence}.}

\section{The labeled case}\label{Sec_L}

Let $H_n$ denote the number of possible
hierarchical orderings or societies, with
exponential generating function (or e.g.f.)
$H(x) = \sum_{n \ge 0} H_n x^n/n!$.
An explicit formula for $H(x)$ follows from
a standard application of the exponential formula
in combinatorics.

\begin{theorem}\label{Th1}
\beql{Eq1}
H(x) ~=~ \exp(B(x) - 1)~,
\eeq
where
$B(x) = 1/(2 - e^x)$ is the e.g.f. for the ordered Bell numbers.
\end{theorem}

\paragraph{Proof.}
This is immediate from (for example) \cite[Cor. 3.4.1]{Wilf94},
or \cite[Chap.~1.4, p.~46]{Berg98}.~~~$\bsq$

The first few values $H_n$ for $n = 0, 1, 2, \ldots$ are
$1, 1, 4, 23, 173, 1602, 17575, 222497, 3188806,$
$50988405, \ldots$ (this is now sequence \htmladdnormallink{A75729}{http://www.research.att.com/cgi-bin/access.cgi/as/njas/sequences/eisA.cgi?Anum=A075729} in \cite{OEIS}).
Table I illustrates the case $n = 3$.

\small \normalsize
\begin{table}\label{TABLEI}
\begin{center}
\caption{For $n = 3$ elements there are $23$ hierarchical orderings
into at most three different subsets.  The subsets are separated by bars,
and the hierarchy within a subset is indicated by the
vertical arrangement.}
\end{center}
\vspace{0.3cm}
\begin{center}
\begin{tabular}{ccc}  \label{n=3}
%
%
\fbox{\begin{tabular*}{50.0pt}{cc}   & \\  & \\ \hspace{0.25cm} 123 \end{tabular*}} & & \\
& & \\
\fbox{\begin{tabular*}{50.0pt}{cc}   & \\ \hspace{0.25cm} 1 & \\ \hspace{0.25cm} 23 &  \end{tabular*}} &
\fbox{\begin{tabular*}{50.0pt}{cc}   & \\ \hspace{0.25cm} 3 & \\ \hspace{0.25cm} 12 &  \end{tabular*}} &
\fbox{\begin{tabular*}{50.0pt}{cc}   & \\ \hspace{0.25cm} 2 & \\ \hspace{0.25cm} 31 &  \end{tabular*}} \\
& & \\
\fbox{\begin{tabular*}{50.0pt}{cc}   & \\  & 31 \\ & 2 \end{tabular*}} &
\fbox{\begin{tabular*}{50.0pt}{cc}   & \\  & 23 \\ & 1 \end{tabular*}} &
\fbox{\begin{tabular*}{50.0pt}{cc}   & \\  & 12 \\ & 3 \end{tabular*}} \\
& & \\
\fbox{\begin{tabular*}{50.0pt}{cc} \hspace{0.25cm} 1 & \\ \hspace{0.25cm} 2 & \\ \hspace{0.25cm} 3 & \end{tabular*}} &
\fbox{\begin{tabular*}{50.0pt}{cc} \hspace{0.25cm} 3 & \\ \hspace{0.25cm} 1 & \\ \hspace{0.25cm} 2 & \end{tabular*}} &
\fbox{\begin{tabular*}{50.0pt}{cc} \hspace{0.25cm} 2 & \\ \hspace{0.25cm} 3 & \\ \hspace{0.25cm} 1 & \end{tabular*}} \\
& & \\
\fbox{\begin{tabular*}{50.0pt}{cc} \hspace{0.25cm} 3 & \\ \hspace{0.25cm} 2 & \\ \hspace{0.25cm} 1 & \end{tabular*}} &
\fbox{\begin{tabular*}{50.0pt}{cc} \hspace{0.25cm} 2 & \\ \hspace{0.25cm} 1 & \\ \hspace{0.25cm} 3 & \end{tabular*}} &
\fbox{\begin{tabular*}{50.0pt}{cc} \hspace{0.25cm} 1 & \\ \hspace{0.25cm} 3 & \\ \hspace{0.25cm} 2 & \end{tabular*}} \\
& & \\
%
%
\fbox{\begin{tabular*}{50.0pt}{c|c}    & \\  &  \\  12 & \hspace{0.05cm} 3 \end{tabular*}} &
\fbox{\begin{tabular*}{50.0pt}{c|c}    & \\  &  \\  31 & \hspace{0.05cm} 2 \end{tabular*}} &
\fbox{\begin{tabular*}{50.0pt}{c|c}    & \\  &  \\  23 & \hspace{0.05cm} 1 \end{tabular*}} \\
& & \\
%
\fbox{\begin{tabular*}{50.0pt}{c|c}    & \\  1 & \\ 2 & \hspace{0.075cm} 3 \end{tabular*}} &
\fbox{\begin{tabular*}{50.0pt}{c|c}    & \\  3 & \\ 1 & \hspace{0.075cm} 2 \end{tabular*}} &
\fbox{\begin{tabular*}{50.0pt}{c|c}    & \\  2 & \\ 3 & \hspace{0.075cm} 1 \end{tabular*}} \\
& & \\
\fbox{\begin{tabular*}{50.0pt}{c|c}   & \\ 3 & \\ 2 & \hspace{0.075cm} 1 \end{tabular*}} &
\fbox{\begin{tabular*}{50.0pt}{c|c}   & \\ 1 & \\ 3 & \hspace{0.075cm} 2 \end{tabular*}} &
\fbox{\begin{tabular*}{50.0pt}{c|c}   & \\ 2 & \\ 1 & \hspace{0.075cm} 3 \end{tabular*}} \\
& & \\
%
%
\fbox{\begin{tabular*}{50.0pt}{cc}  \hspace{-0.19cm} \vline \hspace{0.425cm}  \vline & \\  \hspace{-0.19cm} \vline \hspace{0.425cm} \vline & \\  1 \vline \hspace{0.05cm} 2 \vline \hspace{0.05cm} 3 \end{tabular*}}
\end{tabular}
\end{center}
\end{table}
\renewcommand{\baselinestretch}{2.0} 

Several properties can be deduced from Theorem \ref{Th1}.

(i) By logarithmic differentiation of \eqn{Eq1} 
(cf. \cite[Chap.~1, p.~22]{Wilf94}) we
obtain a recurrence
\beql{rec}
H_n  = \sum_{k = 1}^{n} {n -1 \choose k - 1} B_k H_{n-k} ~. 
\eeq

(ii) Expanding the right-hand
side of \eqn{Eq1} leads to an explicit formula:
\beql{Count}
H_n = \sum_{ (m_1, m_2, \ldots, m_n)}
\frac{ n! \; \prod_{j=1}^{n} B_j^{m_{j}} }{ \prod_{j=1}^{n} m_{j}! \; (j!)^{m_{j}}}   ~,
\eeq
where the sum is over all $(m_1, m_2, \ldots, m_n)$ such
that $\sum_{j=1}^{n} j m_j = n$; that is,
$m_j$ is the number
of subsets of $S$ containing $j$ elements, for $j=1, \ldots, n$.

When $n = 6$, for example,

\begin{math}
H_6 ~=~ \frac{6!}{6! (1!)^6 } \; 1 + 
\frac{6!}{4! 1! (1!)^4 (2!)^1 } \; 1^4 \; 3^1 + 
\frac{6!}{2! 2! (1!)^2 (2!)^2 } \; 1^2 \; 3^2 + 
\frac{6!}{3! (2!)^3 } \; 3^3 + 
\frac{6!}{3! 1! (1!)^3 (3!)^1 } \; 1^3 \; 13^1 + 
\frac{6!}{1! 1! 1! (1!)^1 (2!)^1 (3!)^1 } \; 1^1 \; 3^1 \; 13^1 + 
\frac{6!}{2! 1! (3!)^2 } \; 13^2 + 
\frac{6!}{2! 1! (1!)^2 (4!)^1 } \; 1^2 \; 75^1 + 
\frac{6!}{1! 1! (2!)^1 (4!)^1 } \; 3^1 \; 75^1 + 
\frac{6!}{1! 1! (1!)^1 (5!)^1 } \; 1^1 \; 541^1 + 
\frac{6!}{1! (6!)^1 } \; 4683^1 = 17575.
\end{math}

(iii) The e.g.f. satisfies the differential equation
$$
H'(x) \frac{(2  - \exp(x))^2}{\exp(x)} ~ = ~ H(x) ~.
$$

(iv) Asymptotic behavior. The e.g.f. has an essential singularity at $x = \log 2$.
The asymptotics are sufficiently complicated that we sought
computer assistance.
Salvy's Maple package {\em gdev} \cite{INRIA}, \cite{Salv89} is specifically
designed for this purpose. The result, after
simplification, is that
\beql{asymp}
H_n  ~\sim~ \frac
{n! \, e^{\sqrt{2n/\log 2}}}
{C^{1/4} \, n^{3/4} \, (\log 2)^n} \mbox{~as~} n \rightarrow \infty ~,
\eeq
where
$$
C ~=~ 32 \,  \pi^2 \, \exp(3-\frac{1}{\log 2}) \, \log 2 ~=~ 1038.97\ldots   ~.
$$
This implies that
$$
\log H_n ~\sim~ n\log n -n(1+\log\log 2) + \sqrt{\frac{2n}{\log 2}} + O(\log n) ~.
$$
For comparison, $B_n$ satisfies
\beql{asympB}
B_n  ~\sim~ \frac
{n!}{2 \, (\log 2)^{n+1}} ~,
\eeq
$$
\log B_n ~\sim~ n\log n -n(1+\log\log 2) + O(\log n) ~.
$$

By iterating Theorem \ref{Th1} we can also count
hierarchies of hierarchical orderings. The e.g.f. is 
\beql{HHO}
\exp\left(  \exp\left( \frac{1}{2 - \exp(x)} - 1 \right) - 1 \right) ~,
\eeq
and the first few terms are
$ 1, 1, 6, 52, 588, 8174, 134537, 2554647, 54909468, 1316675221, \ldots$
(sequence \htmladdnormallink{A75756}{http://www.research.att.com/cgi-bin/access.cgi/as/njas/sequences/eisA.cgi?Anum=A075756} in \cite{OEIS}).

\section{The unlabeled case}\label{Sec_U}

If the initial $n$ elements are unlabeled, ordered set-partitions
are called ``compositions'' of $n$, and their number is
$2^{n-1}$ \cite[p. 124]{Rior58}. 
To obtain a hierarchical ordering we partition the
elements into $m \le n$ unlabeled and nonempty subsets,
and form a composition of each subset.
Let $U_n$ denote the number of 
such hierarchical orderings, with
ordinary generating function (or o.g.f.)
$U(x) = \sum_{n \ge 0} U_n x^n$.

\begin{theorem}\label{Th2} 
\beql{EqU1}
U(x) ~=~ \prod_{j \ge 1} \frac{1}{(1-x^j)^{2^{j-1}}} ~.
\eeq
\end{theorem}

\paragraph{Proof.}
This is the unlabeled analogue of Theorem \ref{Th1}.
If $a_n$, $n \ge 1$, is the number of $n$-element
objects with a certain property, then $b_n$,
the number of disjoint unions of such objects
with a total of $n$ elements, where
the order of the components is unimportant,
is given by
$$
1 + \sum_{n=1}^{\infty} b_n x^n ~=~ \prod_{j=1}^{\infty} (1-x^j)^{-a_j}
$$
(see for example \cite[p. 91]{Came89}).~~~$\bsq$

The first few values $U_n$ are
$1, 1, 3, 7, 18, 42, 104, 244, 585, 1373,  \ldots$
(sequence \htmladdnormallink{A34691}{http://www.research.att.com/cgi-bin/access.cgi/as/njas/sequences/eisA.cgi?Anum=A034691} in \cite{OEIS}).
When $n=3$, for example,
the $H_3 = 23$ hierarchical orderings in 
Table I reduce to $U_3 = 7$ when the labels are removed.

{\bf Properties.}
(i) Logarithmic differentiation of \eqn{EqU1} leads to a 
recurrence:
\beql{EqU2}
U_n ~=~ \frac{1}{n} \sum_{k=1}^n \alpha_k U_{n-k} \,,
\mbox{~where~}
\alpha_k ~=~ \sum_{d | k} d 2^{d-1} ~.
\eeq

(ii) It is not so straightforward to find
the asymptotic behavior in this case.  This is to be expected, since
the generating function for the number of 
partitions of $n$ has a similar form to \eqn{EqU1}.
Also, \eqn{EqU1} does not belong to the family
of generalized partition functions considered by Meinardus
and discussed in \cite[Chap. 7]{Andr76}.
However, the saddle point method applies.

\begin{theorem}\label{Th3}
\beql{EqU3}
U_n  ~\sim~ \frac
{2^n e^{\sqrt{2n}}}
{\sqrt{2 \pi} \, 2^{3/4} \, e^{1/4} \, n^{3/4} } 
 \mbox{~as~} n \rightarrow \infty ~.
\eeq
\end{theorem}

\paragraph{Proof.}
We have
\begin{eqnarray*}
\log U(x) & ~=~ & \sum_{n=1}^{\infty} \sum_{m=1}^{\infty} 
\frac{2^{n-1} \, x^{mn} }{m} \\
 & ~=~ & \sum_{N=1}^{\infty} x^N \sum_{ d | N }
\frac{d \, 2^{d-1}}{N} ~,
\end{eqnarray*}
the interchange of summations being justified since
all terms are positive, so
\beql{EqU4}
\log U(x) ~=~ \sum_{k = 1}^{\infty} \frac{x^k}{1-2x^k} ~.
\eeq
This has poles at $2^{-1}, 2^{-1/2}, 2^{-1/3}, \ldots,$
and the radius of convergence is $1/2$.
The pole at $x = 1/2$ dominates, and we apply
the saddle-point method as in \cite[\S 12]{Odly95}.
The saddle point is at $x = r_n$, the solution to
$$
r_n \, \frac{U'(r_n)}{U(r_n)} ~=~ n ~,
$$
which is
$$
r_n ~=~ \frac{1}{2} - \frac{1}{8n} \sqrt{8n+1} + \frac{1}{8n} + O(n^{-3/2}) ~.
$$
Then \cite[Eq. (12.9)]{Odly95} leads to \eqn{EqU3}.~~~$\bsq$

\section{The structure of a random hierarchy}\label{Sec_H}
In this section we consider the case of a single (labeled
or unlabeled) hierarchy.
The elements at the bottom of the hierarchy
will be said to have rank $1$, those at the next level rank $2$, 
and so on.  The maximal rank in a hierarchy is its {\em height}.

Suppose there are $n$ labeled elements. Let $X$ be a
hierarchy chosen at random from the $B_n$ possibilities,
and let $x \in X$ be a randomly chosen element.
There are $h! \left \{  {n \atop h} \right \}$ ways that
$X$ can have height $h$,
where $\left \{  {n \atop h} \right \}$
is a Stirling number of the second kind (cf. \cite[Chap. 7, Problem 44]{GKP94}),
and indeed
$$
B_n ~=~ \sum_{h=0}^{n} h! \left \{  {n \atop h} \right \} ~.
$$
Given  that $X$ has height $h$, $x \in X$ is equally likely
to have any rank from $1$ to $r$.
It follows that the probability that a randomly chosen $x$ has rank $r$
is 
\beql{Eqt7}
P_{n,r}  ~=~ \sum_{i=r}^{n} \frac{1}{i} \,
\frac{i! \, \left \{  {n \atop i} \right \}}{B_n} 
 ~=~  \frac{1}{B_n} \sum_{i=r}^{n} (i-1)! \left \{  {n \atop i} \right \} ~.
\eeq
The average rank is
\beql{Eqt1}
a_n ~=~ \sum_{r=1}^{n} r P_{n,r} ~=~ 
\frac{1}{2 B_n} \, \sum_{i=1}^n (i+1)! \left \{  {n \atop i} \right \}
~.
\eeq 
The numbers $a_n B_n$, $n=0, 1, \ldots$, have e.g.f.
$$
\sum_{n=0}^{\infty} \frac{x^n}{n!} \,
 \sum_{i=1}^{n} \frac{(i+1)!}{2}  \left \{  {n \atop i} \right \}
~=~
- \frac{1}{2} 
\frac{ (e^x-1)(e^x-3)}{(e^x-2)^2} ~,
$$
after some simplification, using standard properties
of Stirling numbers (cf. \cite[Table 351]{GKP94}).
This implies
\beql{Eqt2}
~ a_n B_n   \sim  n! \, \frac{n}{8 (\log 2)^{n+2}} \mbox{~as~} n \rightarrow \infty ~.
\eeq
Combining this with \eqn{asympB} we find that the
average rank is asymptotic to
\beql{Eqt3}
 \frac{n}{4 \log 2} ~=~ 0.36067\ldots n  \mbox{~as~} n \rightarrow \infty ~.
\eeq

In the unlabeled case, if one of the $2^{n-1}$ compositions of $n$
is chosen at random, $X$ (say), and $x \in X$ is picked
at random, one can show that the probability that $x$ has
rank $r$ is
\beql{Eqt4}
\frac{1}{n \, 2^{n-1}} \, \sum_{i=r}^{n}
\left (  {n \atop i} \right )  ~,
\eeq
and the average rank is
$(n+3)/4$.
Ranks are higher in the labeled case.

\paragraph{Acknowledgements}
We thank Jeffrey Lagarias and Andrew Odlyzko
for helpful discussions.
T.W. thanks the Mathematics Department of the University Kassel for
providing access to their computer pool.

\end{document}